\documentclass[11pt,a4paper]{article}
\usepackage{amsmath,amssymb,amsfonts}%
\usepackage{amsthm}%
\usepackage{mathabx}%
\usepackage[title]{appendix}%
\usepackage{xcolor}%
\usepackage{hyperref}%
\usepackage{a4wide}%

\usepackage{authblk}

\title{A short note on inf-sup conditions for the\\ Taylor-Hood family $Q_k$--$Q_{k-1}$}

\author{Walter Zulehner}

\affil{Johann Radon Institute for Computational and Applied Mathematics,\\
Austrian Academy of Sciences, Altenberger Stra{\ss}e 69, 4040 Linz, Austria.\\
E-mail address: \texttt{zulehner@numa.uni-linz.ac.at}}

\date{}

\newtheorem{lem}{Lemma}
\newtheorem{cor}{Corollary}
\newtheorem{thm}{Theorem}

\newcommand{\grad}{\nabla}
\renewcommand{\div}{\operatorname{div}}
\newcommand{\cof}{\operatorname{cof}}

\newcommand{\node}{{\boldsymbol a}}
\newcommand{\noderef}{{\hat{\boldsymbol a}}}
\newcommand{\xref}{{\hat x}}
\newcommand{\sref}{{\hat \sigma}}

\allowdisplaybreaks

\begin{document}

\maketitle

\abstract{We discuss two types of discrete inf-sup conditions of the Stokes problem for the Taylor-Hood family $Q_k$--$Q_{k-1}$ for all $k\in \mathbb{N}$ with $k\ge 2$ in 2D and 3D. While in 2D all results hold for a general class of quadrilateral meshes, the results in 3D are restricted to meshes of parallelepipeds. The analysis is based on an element-wise technique as opposed to the widely used macroelement technique.
It covers Stokes problems with pure Dirichlet boundary conditions, pure Neumann boundary conditions, and mixed boundary conditions.
}

\section{Introduction}\label{sec1}

Let $\Omega \subset \mathbb{R}^d$ , $d \in \{2,3\}$, be a bounded, connected and open set with Lipschitz continuous boundary $\partial \Omega$.
We consider the following Stokes problem: Find the velocity $u$ and the pressure $p$ such that
\begin{alignat*}{2}
  -\mu \, \Delta u + \grad p = f,  \quad 
                      \div u & = 0 & \quad & \text{in} \ \Omega , \\
                      u & = 0 & \quad & \text{on} \ \Gamma_D, \\
   \mu \, \partial_n u - p \, n & = 0 & \quad & \text{on} \ \Gamma_N ,
\end{alignat*}
where $\mu$ is the dynamic viscosity, $f$ is a given force density, $\partial \Omega = \Gamma_D \cup \Gamma_N$ with $\Gamma_D \cap \Gamma_N = \emptyset$, $\partial_n$ indicates the directional derivative with respect to the unit outer normal vector  $n$ to $\partial \Omega$. 
The boundary conditions on $\Gamma_N$ are also known as "do-nothing" boundary conditions, see \cite{Gresho_1991}, \cite{HeRaTu96}.

The standard mixed variational formulation of this problem reads: Find $u \in V$ and $p \in Q$, such that
\begin{alignat*}{3}
    a(u,v) & + b(v,p) &
      & = \langle F,v \rangle &
      \quad & \text{for all} \ v \in V , \\
    b(u,q) & &
      & = 0 & 
      \quad & \text{for all} \ q \in Q,
\end{alignat*}
where
\[
  a(u,v) = \mu 
               \int_\Omega \grad u : \grad v \ d x, \quad
  b(v,q) = - \int_\Omega q \ \div v \ d x,  \quad
  \langle F,v \rangle =\int_\Omega f \cdot v \ d x,
\]
with $f \in L^2(\Omega)^d$ and the function spaces
\begin{enumerate}
\item[(a)]
for pure Dirichlet boundary conditions ($\Gamma_N = \emptyset$):
\[
  V = H_0^1(\Omega)^d, \quad
  Q = \left\{q \in L^2(\Omega) \colon \int_\Omega q \ d x = 0 \right\} =: L_0^2(\Omega),
\]
\item[(b)]
for mixed boundary conditions ($\Gamma_D$, $\Gamma_N$ have positive measure, i.e., $|\Gamma_D|,|\Gamma_N| > 0$):
\[
  V = \{v \in H^1(\Omega)^d \colon v = 0 \ \text{on} \ \Gamma_D \} =: H_{0,\Gamma_D}^1(\Omega)^d, \quad
  Q = L^2(\Omega),
\]
\item[(c)]
for pure  Neumann boundary conditions ($\Gamma_D = \emptyset$):
\[
  V = \left\{v \in H^1(\Omega)^d \colon \int_\Omega v \ d x  = 0 \right\} =: \hat{H}^1(\Omega)^d, \quad
  Q = L^2(\Omega).
\]
\end{enumerate}
Here and throughout the paper we use $L^2(D)$ and $H^1(D)$ to denote the standard Lebesgue and Sobolev spaces of functions on $D$. Moreover, $H_0^1(D)$ is the subspace of functions in $H^1(D)$ with vanishing trace on $\partial D$, its dual space is indicated by $H^{-1}(D)$. The inner product in $L^2(D)$ is denoted by $(.,.)_{0,D}$, the norms in $L^2(D)$, $H^1(D)$, $H^{-1}(D)$ are denoted by $\|.\|_{0,D}$, $\|.\|_{1,D}$,$\|.\|_{-1,D}$, respectively. The same symbols are used for inner product and norms 
in the vector-valued versions of the spaces. For $D = \Omega$ we drop the subscript $D$.
The Euclidean vector norm in $\mathbb{R}^d$ is indicated by $\|.\|$ (without any subscript).  

It is well-known that the famous inf-sup condition (also known as LBB condition) is satisfied for this mixed variational problem: There exists a constant $\beta >0$ such that
\begin{equation} \label{classicalinfsup}
  \sup_{0\neq v \in V} \frac{b(v,q)}{\|v\|_1} 
  \ge \beta \,  \|q\|_0 \quad \text{for all} \ q \in Q.
\end{equation}
The LBB condition is essential for the well-posedness of the Stokes problem. Its verification for the case (a) is typically done by an indirect proof using a compactness argument. We are aware of only one reference, namely \cite{Baerland_2017}, which contains a proof for the case $|\Gamma_N| > 0$. In Appendix \ref{secA1} an alternative proof for this case is given, which follows very closely the standard indirect proof.

A popular class of finite element methods for discretizing the Stokes problem is the generalized Taylor-Hood family $P_k$--$P_{k-1}$ on triangular/tetrahedral meshes with continuous approximations for both velocity and pressure. The discrete analog of \eqref{classicalinfsup}, the discrete LBB condition, which is essential for the stability of the discretization and for error estimates, is well understood for this family, see \cite{Boffi94}, \cite{Boffi97}, and also \cite{Boffi13} for additional references. 

In this paper we focus on the related generalized Taylor-Hood family $Q_k$--$Q_{k-1}$ on quadrilateral/hexahedral meshes under the following assumptions.
\begin{enumerate}
\item[(i)]
$\Omega$ is a polygonal/polyhedral domain, i.e., we do not consider domains with curved boundaries. 
\item[(ii)]
$\mathcal{T}_h$ is a subdivision of $\Omega$ into convex quadrilaterals/hexahedrons $K$, which satisfies the standard compatibility conditions. 
\item[(iii)]
Each element $K \in \mathcal{T}_h$ is of the form
\[
  K = F_K(\hat K)
\]
with the reference element $\hat K = [0,1]^d$, bijective and continuously differentiable mappings $F_K   \colon \mathbb{R}^d \to \mathbb{R}^d$ for each $K \in \mathcal{T}_h$, and the following standard shape regularity conditions. There are positive constants $c_i$ such that
\[
  c_1 \, h_K^d \le |\det J_K(\xref)| \le c_2 \, h_K^d
  \quad \text{and} \quad
  c_3 \, h_K^{2} \, I_d \le J_K(\xref)^\top J_K(\xref) \le c_4 \, h_K^{2} \, I_d
\]
for all $\xref \in \hat K$ with the identity matrix $I_d$ in $\mathbb{R}^d$.
Here $h_K$ denotes the diameter of $K$, $J_K$ denotes the Jacobian of $F_K$ with determinant $\det J_K$, for which we assume that $\det J_K(\xref) > 0$. The mesh size of $\mathcal{T}_h$ is introduced by $h = \max \{h_K \colon K \in \mathcal{T}_h \}$.
\item[(iv)]
For each $K \in \mathcal{T}_h$, at least one of its vertices does not lie on $\Gamma_D$. This excludes only some very degenerate meshes.
\end{enumerate}

The family of generalized Taylor-Hood elements is given by the following pairs of finite element spaces for $k \ge 2$:
\begin{equation} \label{TaylorHood}
\begin{alignedat}{1}
  V_h & = \left\{ v \in V \colon v \circ F_K \in Q_k(\hat{K})^d \ \text{for all} \ K \in \mathcal{T}_h \right\}, \\
  Q_h & = \left\{ q \in Q \cap C(\bar{\Omega}) \colon q \circ F_K \in Q_{k-1}(\hat{K})  \ \text{for all} \ K \in \mathcal{T}_h \right\}
\end{alignedat}
\end{equation}
with the following standard notation of polynomial spaces: Let $P_k(D)$
denote the set of polynomials on $D$ with total degree $\le k$,
and let $P_{k_1,\ldots,k_d}(D)$ denote the set of polynomials on $D$ with degree $\le k_i$ with respect to the $i^\text{th}$ variable, for $i=1,\ldots,d$. Finally, we set $Q_k(D) = P_{k,\ldots,k}(D)$.

The discrete LBB condition associated to these finite element spaces reads as follows: 
\begin{equation} \label{LBBinfsup}
  \sup_{0\neq v_h \in V_h} \frac{b(v_h,q_h)}{\|v_h\|_1} 
  \gtrsim \|q_h\|_0 \quad \text{for all} \ q_h \in Q_h.
\end{equation}
Here and throughout the paper we use the following notation: $A \lesssim B$ iff there is a constant $c > 0$ such that $A \le c \, B$, $A \gtrsim B$ iff $B \lesssim A$, and $A \sim B$ iff $A \lesssim B$ and $A \gtrsim B$. The involved constants are meant to be independent of the mesh size $h$ (or the element size $h_K$ for local estimates on some element $K$), but may depend on $k$ and the shape regularity constants.

A natural condition for the local mappings $F_K$ is
\begin{equation} \label{restrictionQ1}
  F_K \in Q_1(\hat K)^d \quad \text{for all} \ K \in \mathcal{T}_h, 
\end{equation}
in order to allow for general quadrilateral/hexahedral meshes on polygonal/polyhedral domains.
Under Condition \eqref{restrictionQ1} the discrete LBB condition was shown for $d=2$ in \cite{Stenberg_1990}.
We are not aware of any valid proof of the discrete LBB condition for the Taylor-Hood elements for $d=3$ under  Condition \eqref{restrictionQ1}. However, the proof of the discrete LBB condition in \cite{Stenberg_1990} extends to $d=3$ under the stronger condition $F_K \in P_1(\hat K)^3$ for all $K \in \mathcal{T}_h$.
A consequence of this stronger restriction is that for $d=3$ only meshes of parallelepipeds are covered.
The discrete LBB condition could also be shown for the isogeometric generalized Taylor-Hood family, see \cite{Bressan_2011}, \cite{Bressan_2013}. 
The proof there relies on a continuously differentiable parametrization of the domain $\Omega$ on each of a fixed number of patches, which does not cover general quadrilateral/hexahedral meshes.

In this paper we will stick to the same conditions as needed in \cite{Stenberg_1990}, namely
\begin{equation} \label{restriction}
  F_K 
    \in \begin{cases}
          Q_1(\hat K)^2 & \quad \text{for} \ d = 2,\\
          P_1(\hat K)^3 & \quad \text{for} \ d = 3
        \end{cases}
\end{equation}
for all $K \in \mathcal{T}_h$,
since we rely as in \cite{Stenberg_1990} on the exact representation of the bilinear form $b(v,q)$ on $V_h \times Q_h$ by the composite Gauss-Lobatto quadrature rule. In more details, we have
\[
  b(v,q) = \sum_{K\in\mathcal{T}_h} b_K(v,q)
\]
with
\begin{align}
  b_K(v,q)
   & = -\int_K q \, \div v \ d x 
     = \int_K v \cdot \grad q \ d x - \int_{\partial K} v \cdot n \, q \ d \sigma \nonumber \\
   & = \int_{\hat{K}} \hat{v}_K \cdot J_K^{-\top} \grad \hat{q}_K \, \det J_K \ d \xref
        - \int_{\partial \hat K} \hat{v}_K \cdot J_K^{-\top} \hat{n}_K \, \hat{q}_K \, \det J_K \  d \hat\sigma
       \nonumber \\
   & = \int_{\hat{K}} \hat{v}_K \cdot \cof(J_K) \grad \hat{q}_K \ d \xref
        - \int_{\partial \hat K} \hat{v}_K \cdot \cof(J_K) \hat{n}_K \, \hat{q}_K \ d \sref
       \label{integrals}
\end{align}
for all $v \in H^1(K)^d$ and $q \in H^1(K)$
with $\hat{v}_K = v \circ F_K$, $\hat{q}_K = q \circ F_K$, and the outward unit normal vector $\hat{n}_K$ to $\partial \hat K$. 
Here we use the $(\det J) \, J^{-\top} = \cof(J)$, where $\cof(J)$ denotes the cofactor matrix of a matrix $J$.
By elementary calculations it can be shown that Condition \eqref{restriction} implies that 
\[
  \hat{v}_K \cdot \cof(J_K) \grad \hat{q}_K 
    \in Q_{2k-1}(\hat K), \quad
  \hat{v}_K \cdot \cof(J_K) \hat{n}_K \, \hat{q}_K \in Q_{2k-1}(\hat f)
\]
for all $\hat{v}_K \in Q_{k}(\hat K)^d$, $\hat{q}_K \in Q_{k-1}(\hat K)$, and on each edge/face $\hat f \subset \partial \hat K$. These conditions allow to represent the integrals in \eqref{integrals} 
for $\hat{v}_K \in Q_k(\hat K)^d$ and $\hat{q}_K \in Q_{k-1}(\hat K)$ 
by the tensor product $(k+1)$-point Gauss-Lobatto quadrature rule exactly. In the case $d=3$
the weaker Condition \eqref{restrictionQ1} only implies that $\hat{v}_K \cdot \cof(J_K) \grad \hat{q}_K \in Q_{2k}(\hat K)$, $\hat{v}_K \cdot \cof(J_K) \hat{n}_K \, \hat{q}_K \in Q_{2k}(\hat f)$ and one loses the exact representation by the tensor product $(k+1)$-point Gauss-Lobatto quadrature, in general.

In this paper we present an alternative proof for the discrete LBB condition, which covers not only the case of pure Dirichlet conditions as it was previously done but also mixed boundary conditions and pure Neumann boundary conditions.

Moreover, we will also study a second type of discrete inf-sup condition: 
\begin{equation} \label{BPinfsup}
  \sup_{0\neq v_h \in V_h} \frac{b(v_h,q_h)}{\|v_h\|_0} 
  \gtrsim \|\grad q_h\|_0 \quad \text{for all} \ q_h \in Q_h.  
\end{equation}
This condition was discussed by Bercovier and Pironneau in \cite{BercovierPironneau_1979}, which turned out as an enabler of \eqref{LBBinfsup}, see \cite{Verfuerth_1984}. We will refer to this inf-sup condition as the discrete BP condition.
In \cite{BercovierPironneau_1979} the proof was given for $k=2$ and for meshes made of rectangles for $d=2$  and bricks for $d=3$. We will extend this result for all $k \ge 2$ and all meshes satisfying Condition \eqref{restriction}.
Our interest in \eqref{BPinfsup} is not its role for showing the discrete LBB condition but its role in the analysis of the smoothing property of some multigrid method for the Stokes problem, see \cite{jodlbauer2022matrixfree}.

The rest of the paper is organized as follows. In Section \ref{sec2} the technique of $T$-coercivity is discussed, which provides important auxiliary results for Section \ref{sec3}, which is the main section of the paper and contains the analysis of the discrete inf-sup conditions. In Appendix \ref{secA1} known results on the continuous LBB condition are recalled and commented. Appendix \ref{secA2} contains helpful relations used throughout the paper.

From now on the subscript $K$ in $F_K$, $J_K$, $\hat{v}_K$, $\hat{q}_K$, $\hat{n}_K$ is frequently omitted for simplicity, if is clear from the context, which $K$ is meant.

\section{$T$-coercivity} \label{sec2}

A standard technique to show a discrete inf-sup condition is the macroelement technique, see \cite{Stenberg_1990}.
A macroelement $M$ is a connected set of (a few) elements $K \in \mathcal{T}_h$. An essential step of this technique is the verification of an inf-sup condition for pairs  $V_M$, $Q_M$ of local finite element spaces associated to a collection of macroelements $M$ which cover $\Omega$.

In this paper we will consider instead the pairs
\begin{align*}
  \bar{V}_K & = \left\{ v \in H^1(K)^d \colon v \circ F_K \in Q_k(\hat{K})^d, \ v = 0 \ \text{on} \ \partial K \cap \Gamma_D \right\}, \\
  \bar{Q}_K & = \left\{ q \in C(K) \colon q \circ F_K \in Q_{k-1}(\hat{K}) \right\}
\end{align*}
associated to each element $K \in \mathcal{T}_h$ as local counterparts of the spaces
\begin{align*}
  \bar{V}_h & = \left\{ v \in H_{0,\Gamma_D}^1(\Omega)^d \colon v \circ F_K \in Q_k(\hat{K})^d \ \text{for all} \ K \in \mathcal{T}_h \right\}, \\
  \bar{Q}_h & = \left\{ q \in C(\bar{\Omega}) \colon q \circ F_K \in Q_{k-1}(\hat{K})  \ \text{for all} \ K \in \mathcal{T}_h \right\} .
\end{align*}
Note that the spaces $\bar{V}_h$, $\bar{Q}_h$ coincide with the spaces $V_h$, $Q_h$ up to the scaling conditions for pure Dirichlet and pure Neumann boundary conditions. 

The following analysis was stimulated by the discussion in \cite{Wieners_2003}, where Taylor-Hood elements for $k=2$ on meshes of mixed type consisting of hexahedrons, prisms, pyramids, and tetrahedrons are considered, and where it is claimed that, for all $K \in \mathcal{T}_h$, there is a linear mapping $T_K \colon \bar{Q}_K \to \bar{V}_K$ such that
\[ 
  b_K(T_K q, q) \gtrsim h_K^2 \, \|\grad q\|_{0,K}^2
  \quad \text{and} \quad
  \|T_K q\|_{1,K} \lesssim h_K \, \|\grad q\|_{0,K} 
\]
for all $q \in \bar{Q}_K$, see {\cite[Lemma 3]{Wieners_2003}}.
In literature such a property is often called $T$-coercivity.
The proof was given  for $k=2$ and only on reference elements $\hat{K}$. 

Here this result is extended to all $k \ge 2$ and all $K \in \mathcal{T}_h$ but only for quadrilateral/hexahedral elements.

\begin{lem} \label{Wienerslocal}
{Let $k \in \mathbb{N}$ with $k \ge 2$, let $K \in \mathcal{T}_h$ and assume that Condition \eqref{restriction} holds for $K$.}
Then there is a linear mapping $T_K \colon \bar{Q}_K \to \bar{V}_K$ such that
\[
  b_K(T_K q, q) \gtrsim h_K^2 \, \|\grad q\|_{0,K}^2
\]
and
\[
  \|T_K q\|_{0,K} \lesssim h_K^2 \, \|\grad q\|_{0,K},
  \quad
  \|T_K q\|_{1,K} \lesssim h_K \, \|\grad q\|_{0,K}
\]
for all $q \in \bar{Q}_K$.
\end{lem}

\begin{proof}
We will restrict ourselves to the case $d=3$, the case $d=2$ is a simple adaptation of the arguments used for $d=3$ and is omitted.

First of all, we introduce a series of notation related to the local degrees of freedom of $\hat{v} \in Q_k(\hat{K})^d$ and $v \in \bar{V}_K$ (for the analysis only, not for the implementation).
\begin{itemize}
\item
Let $\hat{\mathcal{F}}$ be the set of all faces of $\hat{T}$. We have $\hat{\mathcal{F}} = \hat{\mathcal{F}}_{1,2} \cup \hat{\mathcal{F}}_{1,3} \cup \hat{\mathcal{F}}_{2,3}$ with the pairwise disjoint sets $\hat{\mathcal{F}}_{i,j}$ of all faces of $\hat{K}$ which are parallel to the $\xref_i$-axis and the $\xref_j$-axis.
Let $\hat{\mathcal{E}}$ be the set of all edges of $\hat{T}$. We have $\hat{\mathcal{E}} = \hat{\mathcal{E}}_1 \cup \hat{\mathcal{E}}_2 \cup \hat{\mathcal{E}}_3$ with the pairwise disjoint sets $\hat{\mathcal{E}}_{i}$ of all edges of $\hat{K}$ which are parallel to the $\xref_i$-axis.
The set of all vertices of $\hat{K}$ is denoted by $\hat{\mathcal{V}}$.
\item
As nodal degrees of freedom for the set $Q_k(\hat{K})^d$ of shape functions on the reference element $\hat{K}$ we choose the values at the $(k + 1)^3$ nodes of the form $\noderef = (a_{k_1},a_{k_2},a_{k_3})$, $0 \le k_1, k_2, k_3 \le k$, where $0 = a_0 < a_1 < \ldots < a_k = 1$ are the Gauss-Lobatto points in the unit interval.
The set of all such nodes on $\hat K$ is denoted by $\mathcal{N}(\hat K)$.
The sets $\mathcal{N}_\text{int}(\hat{K})$, $\mathcal{N}_\text{int}(\hat{\mathcal{F}}_{i,j})$, and $\mathcal{N}_\text{int}(\hat{\mathcal{E}}_i)$ are introduced as the sets of all nodes in the interior of $\hat{K}$, in the (relative) interior of a face from $\hat{\mathcal{F}}_{i,j}$, and in the (relative) interior of an edge from $\hat{\mathcal{E}}_i$, respectively. Then, obviously, the set $\mathcal{N}(\hat K)$ is the disjoint union of the sets $\mathcal{N}_\text{int}(\hat{K})$, $\mathcal{N}_\text{int}(\hat{\mathcal{F}}_{i,j})$, $\mathcal{N}_\text{int}(\hat{\mathcal{E}}_i)$, and $\mathcal{V}$.
Moreover, let $\hat{f}$ be a face of $\hat{K}$. Then $\mathcal{N}(\hat{f})$ denotes the set of all nodes which lie on $\hat{f}$. 
\item
As nodal degrees of freedom for $\bar{V}_K$ we choose the values at all nodes from $\mathcal{N}(K) = \{\node = F(\noderef) \colon \noderef \in \mathcal{N}(\hat K)\}$. The subsets $\mathcal{N}_\text{int}(K)$, $\mathcal{N}_\text{int}(\mathcal{F}_{i,j})$, and $\mathcal{N}_\text{int}(\mathcal{E}_i)$ of nodes are introduced accordingly as counterparts of the corresponding sets on $\hat{K}$. Note that the sets of nodes $\mathcal{N}_\text{int}(\mathcal{F}_{i,j})$ and $\mathcal{N}_\text{int}(\mathcal{E}_i)$ actually depend on the element $K$, which, however, is suppressed in the notation. 
\item
Finally, to indicate that nodal points on
the boundary $\Gamma_D$ are excluded we write $\mathcal{N}_\text{int}^*(.)$ instead of $\mathcal{N}_\text{int}(.)$.
\end{itemize}

Let $q \in \bar{Q}_K$. The nodal values of the element $v = T_K q \in \bar{V}_K$ are fixed on $K \in \mathcal{T}_h$ as follows.
\[
   v^{(\node)} = 
   \begin{cases}
     \left(J \grad \hat{q}\right)(\noderef)
       & \quad \text{for} \ \node \in \mathcal{N}_\text{int}(K),\\[2ex]
     \left(J_i \, \partial_i \hat{q} + J_j \, \partial_j \hat{q} \right)(\noderef)
       & \quad \text{for} \ \node \in \mathcal{N}_\text{int}^*(\mathcal{F}_{i,j}),
       \\[2ex]
     \left(J_i \, \partial_i \hat{q} \right)(\noderef) 
       & \quad \text{for} \ \node \in \mathcal{N}_\text{int}^*(\mathcal{E}_i),
       \\[2ex]
      0 
       & \quad \text{for all other} \ \node \in \mathcal{N}(K),
    \end{cases}
\]
with $\noderef = F^{-1}(\node)$, $\hat{q} = q \circ F$, the Jacobian $J$ of $F$, whose $\ell^\text{th}$ column is denoted by $J_\ell$.

By using the property $J^\top \cof(J) = (\det J) \,  I_d $, it is easy to see that
\[
  \left(\hat{v} \cdot \cof(J) \grad \hat{q}\right) (\noderef) = 
   \begin{cases}
     \left(\det J \, \|\grad \hat{q}\|^2\right)(\noderef)
       & \quad \text{for} \ \node \in \mathcal{N}_\text{int}(K),\\[2ex]
     \left(\det J \, (|\partial_i \hat{q}|^2 + |\partial_j \hat{q}|^2 )\right)(\noderef)
       & \quad \text{for} \ \node \in \mathcal{N}_\text{int}^*(\mathcal{F}_{i,j}),\\[2ex]
     \left(\det J \, |\partial_i \hat{q}|^2\right)(\noderef)
       & \quad \text{for} \ \node \in \mathcal{N}_\text{int}^*(\mathcal{E}_i),\\[2ex]
     0
       & \quad \text{for all other} \ \node \in \mathcal{N}(K),
   \end{cases}
\]
and 
\[
  \left(\hat v \cdot \cof(J) \hat n \right)(\noderef)
   = 0 \quad \text{for all} \ \noderef \in \partial \hat K,
\]
with $\hat{v} = v \circ F$. 

By using Condition \eqref{restriction}, 
estimates based on the exactness of the Gauss-Lobatto quadrature (see Appendix \ref{secA2} for details),
and $\det J \sim h_K^3$ it follows that
\begin{align}
  & b_K(v,q) 
    = - \int_K q \, \div v \ d x
      =  \int_{\hat{K}} \hat{v} \cdot \cof(J) \grad \hat{q} \ d \xref 
      \sim \sum_{\node \in \mathcal{N}(K)} (\hat{v} \cdot \cof(J) \grad \hat{q})(\noderef) \nonumber \\
  & \quad \sim \sum_{\node \in \mathcal{N}_\text{int}(K)}
        \left(\det J \, \|\grad \hat{q}\|^2\right)(\noderef)  + 
      \sum_{\substack{\node \in \mathcal{N}_\text{int}^*(\mathcal{F}_{i,j})\\i,j \in \{1,2,3\}, i < j}}
        \left(\det J \, (|\partial_i \hat{q}|^2 + |\partial_j \hat{q}|^2 )\right)(\noderef) \nonumber \\
  & \quad \quad \qquad     {} + 
      \sum_{\substack{\node \in \mathcal{N}_\text{int}^*(\mathcal{E}_i)\\i \in \{1,2,3\}   }}
        \left(\det J \, |\partial_i \hat{q}|^2\right)(\noderef) \nonumber \\
  & \quad \sim h_K^3 \, \underbrace{\Big( \sum_{\node \in \mathcal{N}_\text{int}(K)}
        \|\grad \hat{q}(\noderef)\|^2  + 
      \sum_{\substack{\node \in \mathcal{N}_\text{int}^*(\mathcal{F}_{i,j})\\i,j \in \{1,2,3\}, i < j}}
        \left(|\partial_i \hat{q}(\noderef)|^2 + |\partial_j \hat{q}(\noderef)|^2 \right) + 
      \sum_{\substack{\node \in \mathcal{N}_\text{int}^*(\mathcal{E}_i)\\i \in \{1,2,3\}   }}
        |\partial_i \hat{q}(\noderef)|^2 \Big)}_{\displaystyle =: |\hat{q}|_A^2}. \label{seminormA}
\end{align}
For the $L^2$-norm of $\grad q$ on $K$ we obtain
\begin{align}
  \|\grad q\|_{0,K}^2 
    & \sim h_K \, \| \grad \hat{q}\|_{0,\hat K}^2
      \sim h_K \, \underbrace{\sum_{\node \in \mathcal{N}(K)}
       \|\grad \hat{q}(\noderef)\|^2}_{\displaystyle =: |\hat{q}|_B^2}. \label{seminormB}
\end{align}
Note that $|.|_A$ and $|.|_B$ are semi-norms on $Q_{k-1}(\hat K)$.
Next we will study the kernels of these two semi-norms. 
From the assumptions on $K$ it follows that at least one of the vertices of $K$ does not lie on the boundary part $\Gamma_D$. Without loss of generality we assume this vertex corresponds to the node $\noderef = (0,0,0)$.
Then all nodes from the sets 
\begin{align*}
  \mathcal{N}_1 
    & = \{ \noderef =(a_{k_1},a_{k_2},a_{k_3}) \in \mathcal{N}(\hat K) \colon
           1 \le k_1 \le k-1, \ 0 \le k_2 \le k-1, \ 0 \le k_3 \le k-1 \}, \\
  \mathcal{N}_2 
    & = \{ \noderef =(a_{k_1},a_{k_2},a_{k_3}) \in \mathcal{N}(\hat K) \colon
           0 \le k_1 \le k-1, \ 1 \le k_2 \le k-1, \ 0 \le k_3 \le k-1 \}, \\
  \mathcal{N}_3 
    & = \{ \noderef =(a_{k_1},a_{k_2},a_{k_3}) \in \mathcal{N}(\hat K) \colon
           0 \le k_1 \le k-1, \ 0 \le k_2 \le k-1, \ 1 \le k_3 \le k-1 \}
\end{align*}
correspond to nodes in the interior of $K$, or in the interior of some face, or in the interior of some edge of $K$, but none of them corresponds to nodes on the boundary $\Gamma_D$.

Assume now that $|\hat{q}|_A = 0$ for some $\hat{q} \in Q_{k-1}(\hat K)$.
Then
\[
  \partial_1 \hat{q} (\noderef) = 0
  \quad \text{for all} \ \noderef \in \mathcal{N}_1
  \quad \text{and} \quad
  \partial_1 \hat{q} (\noderef) \in P_{k-2,k-1,k-1}(\hat K).
\]
Since $\mathcal{N}_1$ is a contiguous block of $(k-1)\cdot k \cdot k$ nodes, where the function $\partial_1 \hat{q} \in P_{k-1,k,k}(\Hat K)$ vanishes, 
it immediately follows that $\partial_1 \hat{q}$ vanishes everywhere. The same arguments apply to $\partial_2 \hat{q}$ and $\partial_3 \hat{q}$ with the sets $\mathcal{N}_2$ and $\mathcal{N}_3$, respectively, which eventually show that $\grad \hat{q} = 0$, i.e., $\hat{q}$ is constant. 
So the kernel of the semi-norm $|.|_A$ is $Q_0(\hat K)$. It is obvious that the kernel of the semi-norm $|.|_B$ is $Q_0(\hat K)$, too. 
Since both semi-norms are defined on the finite-dimensional space $Q_{k-1}(\hat K)$ and have the same kernel, they are equivalent. Therefore, we immediately obtain from the equivalence of the two semi-norms, \eqref{seminormA}, and \eqref{seminormB} that
\[
  b_K(v,q) 
    \sim h_K^3 \, |\hat{q}|_A^2 
    \sim h_K^3 \, |\hat{q}|_B^2 
    \sim h_K^2 \, \|\grad q\|_{0,K}^2 .
\]
In order to estimate the norms of  $v$ observe that
\[
  \|v^{(\node)}\|^2 \lesssim h_K^2 \, \|\grad \hat{q}(\noderef)\|^2,
\]
which implies
\begin{align*}
  \|v\|_{0,K}^2 
     \sim h_K^3 \, \|\hat{v}\|_{0,\hat K}^2
     & \sim h_K^3 \, \sum_{\node \in \mathcal{N}(K)} \|v^{(\node)}\|^2 \\
     & \lesssim h_K^5 \, \sum_{\node \in \mathcal{N}(K)} \|\grad \hat{q}(\noderef)\|^2 
       = h_K^5 \, |\hat{q}|_B^2
       \sim h_K^4 \, \|\grad q\|_{0,K}^2, \\
  \|v\|_{1,K}^2 
     \lesssim h_K \, \|\hat{v}\|_{1,\hat K}^2
     & \sim h_K \, \sum_{\node \in \mathcal{N}(K)} \|v^{(\node)}\|^2 \\
     & \lesssim h_K^3 \, \sum_{\node \in \mathcal{N}(K)} \|\grad \hat{q}(\noderef)\|^2 
       = h_K^3 \, |\hat{q}|_B^2
       \sim h_K^2 \, \|\grad q\|_{0,K}^2.
\end{align*}
This completes the proof.
\end{proof}

From this local result on $K$ we easily obtain the following global version.

\begin{cor} \label{Wienersglobal}
{Let $k \in \mathbb{N}$ with $k \ge 2$ and assume that Condition \eqref{restriction} holds for all $K \in \mathcal{T}_h$.}
Then there is a linear mapping $T_h \colon \bar{Q}_h \to \bar{V}_h$ such that
\[
  b(T_h q_h, q_h) \gtrsim \sum_{K \in \mathcal{T}_h} h_K^2 \, \|\grad q_h\|_{0,K}^2
\]
and
\[
  \|T_h q_h\|_0   \lesssim \Big(\sum_{K \in \mathcal{T}_h} h_K^4 \, \|\grad q_h\|_{0,K}^2\Big)^{1/2},
  \quad 
  \|T_h q_h\|_1   \lesssim \Big(\sum_{K \in \mathcal{T}_h} h_K^2 \, \|\grad q_h\|_{0,K}^2\Big)^{1/2}
\]
for all $q_h \in \bar{Q}_h$.
\end{cor}

\begin{proof}
Let $q_h \in \bar{Q}_h$.
It is easy to see that the definition of the nodal values of $v$ in the proof of Lemma \ref{Wienerslocal} on common vertices, edges and faces of neighboring elements is consistent. This results in a well-defined and continuous finite element function
\[ 
  v_h = T_h q_h \in \bar{V}_h
  \quad \text{with} \quad v_h|_K = T_K (q_h|_K),
\]
for which we have
\begin{align}
  b(v_h,q_h) & 
   = \sum_{K\in \mathcal{T}_h} b_K(T_K (q_h|_K),q_h) 
   \sim \sum_{K\in \mathcal{T}_h}  h_K^2 \, \|\grad q_h\|_{0,K}^2 
\end{align}
and
\begin{align*}
  \|v_h\|_0^2 
    & = \sum_{K\in\mathcal{T}_h} \|v_h\|_{0,K}^2 
  \lesssim \sum_{K\in\mathcal{T}_h} h_K^4 \, \|\grad q_h\|_{0,K}^2 \\
  \|v_h\|_1^2 
    & = \sum_{K\in\mathcal{T}_h} \|v_h\|_{1,K}^2 
  \lesssim \sum_{K\in\mathcal{T}_h} h_K^2 \, \|\grad q_h\|_{0,K}^2 .
\end{align*}
\end{proof}

The results in Lemma \ref{Wienerslocal} and Corollary \ref{Wienersglobal} suffice to discuss the discrete inf-sup conditions for the case (a) of pure Dirichlet boundary conditions. If $|\Gamma_N| > 0$, we need an  additional result, for which we introduce the following notation
\[
  {\mathcal{T}_h^N} = \{ K \in \mathcal{T}_h \colon f \subset \Gamma_N \ \text{for some face} \ f \ \text{of} \ K \}.
\]
For $K \in {\mathcal{T}_h^N}$, the set of all faces of $K$ which lie on $\Gamma^N$ are denoted by $\mathcal{F}^N(K)$.
The set $\mathcal{N}(\mathcal{F}^N)$ is formed by all nodes (as introduced in the proof of Lemma \ref{Wienerslocal}) on faces from $\mathcal{F}^N(K)$.
The set $\mathcal{N}_\text{int}(\mathcal{F}^N)$ is formed by all nodes in the (relative) interior of faces from $\mathcal{F}^N(K)$.

\begin{lem} \label{WienerslocalGammaN}
{Let $k \in \mathbb{N}$ with $k \ge 2$, let $K \in \mathcal{T}_h^N$ and assume that Condition \eqref{restriction} holds for $K$.}
Then there is a linear mapping $T_K^N \colon \bar{Q}_K \to \bar{V}_K$ such that
\[ 
  b_K(T_K^N q, q) 
  \gtrsim h_K \, \vvvert q\vvvert_{0,\partial K \cap \Gamma_N}^2 
    - h_K^{3/2} \, \vvvert q\vvvert_{0,\partial K \cap \Gamma_N} \, \|\grad q\|_{0,K}
\]
and
\[
  \|T_K^N q\|_{1,K} 
  \lesssim h_K^{1/2} \, \vvvert q\vvvert_{0,\partial K \cap \Gamma_N}
\]
for all $q \in \bar{Q}_K$,
with the ($L^2$-like) semi-norm $\vvvert.\vvvert_{0,\partial K \cap \Gamma_N}$ given by
\[
  \vvvert q\vvvert_{0,\partial K \cap \Gamma_N}^2
    = h_K^{d-1} \, \sum_{\node \in \mathcal{N}_\text{int}(\mathcal{F}^N)}
       |\hat{q}(\noderef)|^2  . 
\]
\end{lem}

\begin{proof}
We will restrict ourselves to the case $d=3$, the case $d=2$ is a simple adaptation of the arguments used for $d=3$ and is omitted.

Let $q \in \bar{Q}_K$. The nodal values of the element $v = T_K^N q \in \bar{V}_K$ are fixed on $K \in \mathcal{T}$ as follows.
\[
   v^{(\node)} = 
   \begin{cases}
      -\left(\hat{q} \, J \hat{n}\right) (\noderef)
       & \quad \text{for} \ \node \in \mathcal{N}_\text{int}(\mathcal{F}^N), \\[1ex]
      0 
       & \quad \text{for all other} \ \node \in \mathcal{N}(K),
    \end{cases}
\]
with $\noderef = F^{-1}(\node)$, $\hat{q} = q \circ F$, the Jacobian $J$ of $F$.

By using the property $J^\top \cof(J) = (\det J) \,  I_d $, it is easy to see that
\[
  \left(\hat{v} \cdot \cof(J) \grad \hat{q}\right) (\noderef) = 
   \begin{cases}
     -\left(\det J \, (\hat{q} \, \hat n \cdot \grad \hat{q})\right)(\noderef)
       & \quad \text{for} \ \node \in \mathcal{N}_\text{int}(\mathcal{F}^N), \\[1ex]
     0
       & \quad \text{for all other} \ \node \in \mathcal{N}(K),
   \end{cases}
\]
and 
\[
  - \left(\hat{v} \cdot \cof(J) \hat n \, q\right)(\noderef)
   = \begin{cases}
     \left(\det J \, |\hat{q}|^2\right) (\noderef)
       & \quad \text{for} \ \node \in \mathcal{N}_\text{int}(\mathcal{F}^N), \\[1ex]
     0
       & \quad \text{for all other} \ \node \in \mathcal{N}(K),
   \end{cases}
\]
with $\hat{v} = v \circ F$.

We have
\begin{align*}
  b_K(v,q) 
    & = \int_{\hat{K}} \hat{v} \cdot \cof(J) \grad \hat{q} \ d \xref 
        - \int_{\partial \hat K \cap \Gamma_N} \hat{v} \cdot \cof(J) \hat n \, \hat{q} \ d \sref .
\end{align*}
From estimates based on the exactness of the Gauss-Lobatto quadrature (see Appendix \ref{secA2} for details),
and $\det J \sim h_K^3$, it follows that
\begin{align}
  - \int_{\partial \hat K \cap \Gamma_N} \hat{v} \cdot \cof(J) \hat n \, \hat{q} \ d \sref 
    & \sim \sum_{\node \in \mathcal{N}_\text{int}(\mathcal{F}^N)} \left(\det J \, |\hat{q}|^2\right) (\noderef) \nonumber \\
    & \sim h_K^3 \, \sum_{\node \in \mathcal{N}_\text{int}(\mathcal{F}^N)}
       |\hat{q}(\noderef)|^2  = h_K \, \vvvert q \vvvert^2 . \label{inequ1}
\end{align}
Since $\hat{v} \cdot \cof(J) \grad \hat{q} \in Q_{2k-1}$, and $\det J \sim h_K^3$, we have
\begin{align}
  \Big| \int_{\hat{K}} \hat{v} \cdot \cof(J) \grad \hat{q} \ d \xref \Big|
   & \lesssim \sum_{\node \in \mathcal{N}_\text{int}(\mathcal{F}^N)} 
       \left(\det J \, |\hat{q}| \, \|\grad \hat{q}\|\right)(\noderef) \nonumber \\
   & \sim h_K^3 \, \sum_{\node \in \mathcal{N}_\text{int}(\mathcal{F}^N)} \left(|\hat{q}|
         \, \|\grad \hat{q}\|\right)(\noderef) \nonumber \\
   & \le h_K^3 \, \Big(\sum_{\node \in \mathcal{N}_\text{int}(\mathcal{F}^N)} |\hat{q}(\noderef)|^2 \Big)^{1/2}
         \Big(\sum_{\node \in \mathcal{N}_\text{int}(\mathcal{F}^N)} \|\grad \hat{q}(\noderef)\|^2\Big)^{1/2} \nonumber \\
   & = h_K^2 \, \vvvert q\vvvert \,
         \Big(\sum_{\node \in \mathcal{N}_\text{int}(\mathcal{F}^N)} \|\grad \hat{q}(\noderef)\|^2\Big)^{1/2}, \label{inequ2}
\end{align}
where the last estimate follows from the Cauchy inequality.
Since
\begin{align}
  \sum_{\node \in \mathcal{N}_\text{int}(\mathcal{F}^N)} \|\grad \hat{q}(\noderef)\|^2
    \le \sum_{\node \in \mathcal{N}(K)} \|\grad \hat{q}(\noderef)\|^2
    = |\hat{q}|_B^2
    \sim h_K^{-1} \, \|\grad q\|_{0,K}^2 , \label{inequ3}
\end{align}
it follows from \eqref{inequ1}, \eqref{inequ2}, \eqref{inequ3} that
\begin{align*}
  & b_K(v,q) 
    \gtrsim h_K \, \vvvert q\vvvert_{0,\partial K \cap \Gamma_N}^2
      - h_K^{3/2} \, \vvvert q\vvvert_{0,\partial K \cap \Gamma_N} \, \|\grad q\|_{0,K}.
\end{align*}
In order to estimate the norm of $v$ observe that
\[
  \|v^{(\node)}\|^2 \lesssim h_K^2 \, |\hat{q}(\noderef)|^2 ,
\]
which implies
\begin{align*}
  \|v\|_{1,K}^2 
    \lesssim h_K \, \|\hat{v}\|_{1,\hat K}^2
    & \lesssim h_K \, \sum_{\node \in \mathcal{N}_\text{int}(\mathcal{F}^N)} \|v^{(\node)}\|^2 \\
    & \lesssim h_K^3 \, \sum_{\node \in \mathcal{N}_\text{int}(\mathcal{F}^N)} |\hat{q}(\noderef)|^2 
      = h_K \, \vvvert q\vvvert_{0,\partial K \cap \Gamma_N}^2.
\end{align*}
This completes the proof.
\end{proof}

From this local result on $K$ we easily obtain the following global version.

\begin{cor} \label{WienersglobalGammaN}
{Let $k \in \mathbb{N}$ with $k \ge 2$ and assume that Condition \eqref{restriction} holds for all $K \in \mathcal{T}_h^N$.}
Then there is a linear mapping $T_h^N \colon \bar{Q}_h \to \bar{V}_h$ such that
\[ 
  b(T_h^N q_h, q_h) 
  \gtrsim \sum_{K \in {\mathcal{T}_h^N}} h_K \, \vvvert q_h\vvvert_{0,\partial K \cap \Gamma_N}^2 
    - \Big(\sum_{K \in {\mathcal{T}_h^N}} h_K \, \vvvert q_h\vvvert_{0,\partial K \cap \Gamma_N}^2\Big)^{1/2} \, 
      \Big(\sum_{K \in {\mathcal{T}_h^N}} h_K^2 \, \|\grad q_h\|_{0,K}^2 \Big)^{1/2}
\]
and
\[
  \|T_h^N q_h\|_1 
  \lesssim \Big(\sum_{K \in {\mathcal{T}_h^N}} h_K \, \vvvert q_h\vvvert_{0,\partial K \cap \Gamma_N}^2\Big)^{1/2}
\]
for all $q_h \in \bar{Q}_h$.
\end{cor}

\begin{proof}
Let $q_h \in \bar{Q}_h$.
It is easy to see that the definition of the nodal values of $v$ in the proof of Lemma \ref{WienerslocalGammaN} on common vertices, edges and faces of neighboring elements is consistent. This results in a well-defined and continuous finite element function
\[ 
  v_h = T_h q_h \in \bar{V}_h
  \quad \text{with} \quad v_h|_K = T_K (q_h|_K),
\]
for which we have
\begin{align*}
  b(v_h,q_h) 
   & = \sum_{K\in {\mathcal{T}_h^N}} b_K(T_K (q_h|_K),q_h) \\
   & \gtrsim \sum_{K\in {\mathcal{T}_h^N}}  h_K \, \vvvert q_h\vvvert_{0,\partial K \cap \Gamma_N}^2 
   - \sum_{K \in {\mathcal{T}_h^N}} h_K^{3/2} \, \vvvert q_h\vvvert_{0,\partial K \cap \Gamma_N} \, \|\grad q_h\|_{0,K} \\
   & \ge \sum_{K \in {\mathcal{T}_h^N}}  h_K \, \vvvert q_h\vvvert_{0,\partial K \cap \Gamma_N}^2 
    - \Big(\sum_{K \in {\mathcal{T}_h^N}}  h_K \, \vvvert q_h\vvvert_{0,\partial K \cap \Gamma_N}^2\Big)^{1/2} \, 
      \Big(\sum_{K \in {\mathcal{T}_h^N}} h_K^2 \, \|\grad q_h\|_{0,K}^2 \Big)^{1/2}
\end{align*}
by using the Cauchy inequality and
\[
  \|v_h\|_1^2 
   = \sum_{K \in {\mathcal{T}_h^N}} \|v_h\|_{1,K}^2 
  \lesssim \sum_{K \in {\mathcal{T}_h^N}} h_K \, \vvvert q_h\vvvert_{0,\partial K \cap \Gamma_N}^2. 
\]
\end{proof}

\section{Discrete inf-sup conditions}\label{sec3}

\subsection{The discrete LBB condition}\label{subsec1}

We start the discussion with an element-wise version of the discrete LBB condition, for which we need the following auxiliary result.

\begin{lem} \label{lemlocal}
{Let $k \in \mathbb{N}$ with $k \ge 2$, let $K \in \mathcal{T}_h$ and assume that Condition \eqref{restriction} holds for $K$.}
Then
\[
  \sup_{0 \neq v \in \bar{V}_K} \frac{ b_K(v,q)}{\|v\|_{1,K}} 
    \gtrsim h_K \, \|\grad q\|_{0,K}
  \quad \text{for all} \ q \in \bar{Q}_K.
\]
If, additionally, $K \in \mathcal{T}_h^N$, then
\[
  \sup_{0 \neq v \in \bar{V}_K}
  \frac{ b_K(v,q)}{\|v\|_{1,K}} 
    \gtrsim \Big(h_K^2 \, \|\grad q\|_{0,K}^2 + h_K \, \|q\|_{0,\partial K \cap \Gamma_N}^2\Big)^{1/2}
  \quad \text{for all} \ q \in \bar{Q}_K.
\]
\end{lem}

\begin{proof}
The first estimate follows immediately from Lemma \ref{Wienerslocal}:
\begin{align}
  \sup_{0 \neq v \in \bar{V}_K} \frac{ b_K(v,q)}{\|v\|_{1,K}}
  \ge \frac{ b_K(T_K q,q)}{\|T_K q\|_{1,K}}
  \gtrsim \frac{h_K^2 \, \|\grad q\|_{0,K}^2}{h_K \, \|\grad q\|_{0,K}}
   = h_K \, \|\grad q\|_{0,K}. \label{firstestlocal}
\end{align}
If $K \in \mathcal{T}_h^N$, then we have by Lemma \ref{WienerslocalGammaN} 
\begin{align}
  \sup_{0 \neq v \in \bar{V}_K} \frac{ b_K(v,q)}{\|v\|_{1,K}}
  \ge \frac{ b_K(T_K^N q,q)}{\|T_K^N q\|_{1,K}}
    & \gtrsim \frac{h_K \vvvert q \vvvert_{0,\partial K \cap \Gamma_N}^2 
    - h_K^{3/2} \, \vvvert q \vvvert_{0,\partial K \cap \Gamma_N} \, \|\grad q\|_{0,K}}{h_K^{1/2} \, \vvvert q \vvvert_{0,\partial K \cap \Gamma_N}} \nonumber \\
    & = h_K^{1/2} \, \vvvert q \vvvert_{0,\partial K \cap \Gamma_N} - h_K \, \|\grad q\|_{0,K}. \label{secondestlocal}
\end{align}
From \eqref{firstestlocal} and \eqref{secondestlocal} we obtain
\begin{align*}
  \sup_{0 \neq v \in \bar{V}_K} \frac{ b_K(v_h,q_h)}{\|v_h\|_{1,K}} 
    & \gtrsim h_K \, \|\grad q\|_{0,K} + h_K^{1/2} \, \vvvert q \vvvert_{0,\partial K \cap \Gamma_N} \\
    & \ge \Big( h_K^2 \, \|\grad q\|_{0,K}^2 + h_K \, \vvvert q \vvvert_{0,\partial K \cap \Gamma_N}^2)\Big)^{1/2}
  \quad \text{for all} \ q \in \bar{Q}_K.
\end{align*}
It remains to show that
\[
  h_K^2 \, \|\grad q\|_{0,K}^2 + h_K \, \vvvert q\vvvert_{0,\partial K \cap \Gamma_N}^2
    \sim h_K^2 \, \|\grad q\|_{0,K}^2 + h_K \, \|q\|_{0,\partial K \cap \Gamma_N}^2 .
\]
Now we have
\begin{align*}
  h_K \, \|q\|_{0,\partial K \cap \Gamma_N}^2
    & = h_K \, \sum_{f \in \mathcal{F}^N(K)} \|q\|_{0,f}^2 
     \sim h_K^3 \, \sum_{f \in \mathcal{F}^N(K)} \|\hat{q}\|_{0,\hat f}^2 
     \sim h_K^3 \, \sum_{\node \in \mathcal{N}(\mathcal{F}^N)} |\hat{q}(\noderef)|^2.
\end{align*}
Together with \eqref{seminormB} we obtain
\begin{align*}
  h_K^2 \, \|\grad q\|_{0,K}^2 + h_K \, \|q\|_{0,\partial K \cap \Gamma_N}^2 
    \sim h_K^3 \, \underbrace{\Big[
     \sum_{\node \in \mathcal{N}(K)} \|\grad \hat{q}(\noderef)\|^2
    + \sum_{\node \in \mathcal{N}(\mathcal{F}^N)} |\hat{q}(\noderef)|^2\Big]}_{\displaystyle =: \|\hat{q}_h\|_C^2}.
\end{align*}
Similarly, we have
\begin{align*}
  h_K^2 \, \|\grad q\|_{0,K}^2 + h_K \, \vvvert q\vvvert_{0,\partial K \cap \Gamma_N}^2 
    \sim h_K^3 \, \underbrace{\Big[  \sum_{\node \in \mathcal{N}(K)} \|\grad \hat{q}(\noderef)\|^2
    + \sum_{\node \in \mathcal{N}_\text{int}(\mathcal{F}^N)} |\hat{q}(\noderef)|^2\Big]}_{\displaystyle =:\|\hat{q}_h\|_D^2}.
\end{align*}
Since both $\|.\|_C$ and $\|.\|_D$ are norms on the finite-dimensional space $Q_{k-1}(\hat K)$, they are equivalent, which completes the proof.
\end{proof}

From Lemma \ref{lemlocal} the element-wise version of the discrete LBB condition easily follows.
\begin{thm}
{Let $k \in \mathbb{N}$ with $k \ge 2$, let $K \in \mathcal{T}_h$ and assume that Condition \eqref{restriction} holds for $K$.}
Then
\[
  \sup_{0 \neq v \in \bar{V}_K} \frac{ b_K(v,q)}{\|v\|_{1,K}} 
    \gtrsim \|q\|_{0,K}
  \quad \text{for all} \ q \in \bar{Q}_K \ \text{with} \ \int_K q \ d x = 0.
\]
If, additionally, $K \in \mathcal{T}_h^N$, then
\[
  \sup_{0 \neq v \in \bar{V}_K}
  \frac{ b_K(v,q)}{\|v\|_{1,K}} 
    \gtrsim  \|q\|_{0,K}
  \quad \text{for all} \ q \in \bar{Q}_K.
\]
\end{thm}

\begin{proof}
The estimates follow immediately from Lemma \ref{lemlocal} and the Friedrichs/Poincar\'{e}-type inequalities 
\begin{align*}
   \|q\|_{0,K}^2 
   \lesssim h_K^{-d} \, \Big( \int_K q(x) \ d x \Big)^2 + h_K^2 \, \|\grad q\|_{0,K}^2
\end{align*}
and
\begin{align*}
   \|q\|_{0,K}^2
     \lesssim h_K \, \|q\|_{0,f}^2 +  h_K^2 \, \|\grad q\|_{0,K}^2,
\end{align*}
where $f$ is some face of $K$ with $f \subset \Gamma_N$.
The proofs of these estimates are analogous to the proof of Lemma 1.1.3 in \cite{Necas12}.
\end{proof}

For the global version of the discrete LBB condition the following global version of Lemma \ref{lemlocal} is needed.

\begin{cor} \label{verf4}
{Let $k \in \mathbb{N}$ with $k \ge 2$ and assume that Condition \eqref{restriction} holds for all $K \in \mathcal{T}_h$.}
Then
\[
  \sup_{0 \neq v_h \in \bar{V}_h} \frac{ b(v_h,q_h)}{\|v_h\|_1} 
    \gtrsim \Big(\sum_{K \in \mathcal{T}_h} h_K^2 \, \|\grad q_h\|_{0,K}^2 + \sum_{K \in {\mathcal{T}_h^N}} h_K \, \|q_h\|_{\partial K \cap \Gamma_N}^2\Big)^{1/2}
  \quad \text{for all} \ q_h \in \bar{Q}_h.
\]
\end{cor}

\begin{proof}
The proof follows closely the proof on Lemma \ref{lemlocal}.
By Corollary \ref{Wienersglobal} we have
\begin{align}
  \sup_{0 \neq v_h \in \bar{V}_h} \frac{ b(v_h,q_h)}{\|v_h\|_1}
    \ge \frac{ b(T_h q_h,q_h)}{\|T_h q_h\|_1}
  \gtrsim \frac{\sum_{K \in \mathcal{T}_h} h_K^2 \, \|\grad q_h\|_{0,K}^2}{\Big(\sum_{K \in \mathcal{T}_h} h_K^2 \, \|\grad q_h\|_{0,K}^2\Big)^{1/2}}
   = \Big(\sum_{K \in \mathcal{T}_h} h_K^2 \, \|\grad q_h\|_{0,K}^2\Big)^{1/2}. \label{firstestglobal}
\end{align}
By Corollary \ref{WienersglobalGammaN} we have
\begin{align}
  & \sup_{0 \neq v_h \in \bar{V}_h} \frac{ b(v_h,q_h)}{\|v_h\|_1}
      \ge \frac{ b(T_h^N q_h,q_h)}{\|T_h^N q_h\|_1} \nonumber \\
  & \quad \gtrsim \frac{\sum_{K \in {\mathcal{T}_h^N}} h_K \, \vvvert q_h\vvvert_{0,\partial K \cap \Gamma_N}^2 
    - \Big(\sum_{K \in {\mathcal{T}_h^N}} h_K \, \vvvert q_h\vvvert_{0,\partial K \cap \Gamma_N}^2\Big)^{1/2} \, 
      \Big(\sum_{K \in {\mathcal{T}_h^N}} h_K^2 \, \|\grad q_h\|_{0,K}^2 \Big)^{1/2}}{\Big(\sum_{K \in {\mathcal{T}_h^N}} h_K \, \vvvert q_h\vvvert_{0,\partial K \cap \Gamma_N}^2\Big)^{1/2}} \nonumber \\
  & \quad = \Big(\sum_{K \in {\mathcal{T}_h^N}} h_K \, \vvvert q_h\vvvert_{0,\partial K \cap \Gamma_N}^2\Big)^{1/2}
            - \Big(\sum_{K \in {\mathcal{T}_h^N}} h_K^2 \, \|\grad q_h\|_{0,K}^2 \Big)^{1/2}. \label{secondestglobal}
\end{align}
From \eqref{firstestglobal} and \eqref{secondestglobal} we obtain
\begin{align*}
  \sup_{0 \neq v_h \in \bar{V}_h} \frac{ b(v_h,q_h)}{\|v_h\|_1} 
   & \gtrsim \Big(\sum_{K \in \mathcal{T}_h} h_K^2 \, \|\grad q_h\|_{0,K}^2\Big)^{1/2}
      + \Big(\sum_{K \in {\mathcal{T}_h^N}} h_K \, \vvvert q_h\vvvert_{0,\partial K \cap \Gamma_N}^2\Big)^{1/2} \\
   & \ge \Big(\sum_{K \in \mathcal{T}_h} h_K^2 \, \|\grad q_h\|_{0,K}^2 
      + \sum_{K \in {\mathcal{T}_h^N}} h_K \, \vvvert q_h\vvvert_{0,\partial K \cap \Gamma_N}^2\Big)^{1/2} \\
   & \sim \Big(\sum_{K \in \mathcal{T}_h} h_K^2 \, \|\grad q_h\|_{0,K}^2 
      + \sum_{K \in {\mathcal{T}_h^N}} h_K \, \| q_h\|_{0,\partial K \cap \Gamma_N}^2\Big)^{1/2}
   \quad \text{for all} \  q_h \in \bar{Q}_h.
\end{align*}
\end{proof}

In order to continue we apply the so-called Verf\"urth trick. We need a slightly modified variant compared to \cite{Verfuerth_1984,Stenberg_1990,Wieners_2003}.

\begin{lem} \label{verf1}
Let $\Omega \subset \mathbb{R}^d$ be a bounded, connected and open set with Lipschitz continuous boundary. Let $V_h \subset V \subset H_{0,\Gamma_D}^1(\Omega)^d$ and $Q_h \subset Q \cap H^1(\Omega)$ with $Q \subset L^2(\Omega)$ be closed subspaces. 
Assume that the LBB condition is satisfied on $V$ and $Q$ and that there exists a linear operator $\Pi_h \colon V \to V_h$
such that
\[
  \Big( \sum_{K \in \mathcal{T}_h}
      h_K^{-2} \, \|v - \Pi_h v\|_{0,K}^2 + 
      \sum_{K \in {\mathcal{T}_h^N}} h_K^{-1} \, \|v-\Pi_h v\|_{0,\partial K \cap \Gamma_N}^2\Big)^{1/2}
   \lesssim \|v\|_1  
\]
and
\[
  \|\Pi_h v\|_1 
  \lesssim \|v\|_1
  \quad \text{for all} \ v \in H^1(\Omega)^d.
\]  
Then
\[
  \sup_{0 \neq v_h \in V_h} \frac{b(v_h,q_h)}{\|v_h\|_1} 
   \gtrsim \| q_h \|_0
    - \Big( \sum_{K \in \mathcal{T}_h} h_K^2 \, \|\grad q_h\|_{0,K}^2 
      + \sum_{K \in {\mathcal{T}_h^N}} h_K \, \|q_h\|_{0,\partial K \cap \Gamma_N}^2\Big)^{1/2}
\]
for all $q_h \in Q_h$.
\end{lem}

\begin{proof}
Let $q_h \in Q_h$. From the continuous LBB condition it follows that there exists a $\bar{v} \in V$ with
\[  
   \frac{b(\bar{v},q_h)}{\|\bar{v}\|_1}
   \gtrsim
   \|q_h\|_0 . 
\]
Then
\begin{align*}
  \sup_{v_h \in V_h} \frac{b(v_h,q_h)}{\|v_h\|_1}
  & \ge \frac{b(\Pi_h \bar{v},q_h)}{\|\Pi_h \bar{v}\|_1} 
    \gtrsim \frac{b(\Pi_h \bar{v},q_h)}{\|\bar{v}\|_1} 
     =  \frac{b(\bar{v},q_h)}{\|\bar{v}\|_1}
	 + \frac{b(\Pi_h \bar{v} - \bar{v},q_h)}%
               {\|\bar{v}\|_1} \\
  & \gtrsim \| q_h \|_0
	 + \frac{b(\Pi_h \bar{v} - \bar{v},q_h)}%
               {\|\bar{v}\|_1} .
\end{align*}
Now we have
\begin{align*}
  & b(\Pi_h \bar{v} - \bar{v},q_h)
    = \int_\Omega \grad q_h \cdot (\Pi_h \bar{v} - \bar{v} ) \ d x
      - \int_{\Gamma_N} q_h \, (\Pi_h \bar{v} - \bar{v}) \cdot n \ d \sigma \\
  & \quad = \sum_{K \in \mathcal{T}_h} \int_K \grad q_h \cdot (\Pi_h \bar{v} - \bar{v} ) \ d x
     - \sum_{K \in {\mathcal{T}_h^N}} \int_{\partial K \cap \Gamma_N} q_h \, (\Pi_h \bar{v} - \bar{v}) \cdot n \ d \sigma \\
  & \quad \ge - \sum_{K \in \mathcal{T}_h} \|\grad q_h\|_{0,K} \, \|\bar{v} - \Pi_h \bar{v}\|_{0,K} 
     - \sum_{K \in {\mathcal{T}_h^N}} \|q_h\|_{0,\partial K \cap \Gamma_N} \, \|\Pi_h \bar{v} - \bar{v}\|_{0,\partial K \cap \Gamma_N} \\
  & \quad \ge - \Big( \sum_{K \in \mathcal{T}_h} h_K^2 \, \|\grad q_h\|_{0,K}^2 
      + \sum_{K \in {\mathcal{T}_h^N}} h_K \, \|q_h\|_{0,\partial K \cap \Gamma_N}^2\Big)^{1/2} \\*
  & \quad \qquad \Big( \sum_{K \in \mathcal{T}_h}
      h_K^{-2} \, \|v - \Pi_h v\|_{0,K}^2 + 
      \sum_{K \in {\mathcal{T}_h^N}} h_K^{-1} \, \|v-\Pi_h v\|_{0,\partial K \cap \Gamma_N}^2\Big)^{1/2} \\
  & \quad \gtrsim - \Big( \sum_{K \in \mathcal{T}_h} h_K^2 \, \|\grad q_h\|_{0,K}^2 
      + \sum_{K \in {\mathcal{T}_h^N}} h_K \, \|q_h\|_{0,\partial K \cap \Gamma_N}^2\Big)^{1/2} \, \|v\|_1 .
\end{align*}
Hence
\begin{align*}
  \sup_{0 \neq v_h \in V_h} \frac{b(v_h,q_h)}{\|v_h\|_1} 
    \gtrsim \| q_h \|_0
	 - 
	 \Big(\sum_{K \in \mathcal{T}_h} h_K^2 \, \|\grad q_h\|_{0,K}^2 + \sum_{K \in {\mathcal{T}_h^N}} h_K \, \|q_h\|_{\partial K \cap \Gamma_N}^2\Big)^{1/2} .
\end{align*}
\end{proof}

Finally, we obtain
\begin{thm} \label{discreteLBB}
{Let $k \in \mathbb{N}$ with $k \ge 2$ and assume that Condition \eqref{restriction} holds for all $K \in \mathcal{T}_h$.}
Then
\[
  \sup_{0 \neq v_h \in V_h}
  \frac{ b(v_h,q_h)}{\|v_h\|_1} \gtrsim \|q_h\|_0 
        \quad \text{for all} \ q_h \in Q_h 
\]
for $V_h$ and $Q_h$ given by \eqref{TaylorHood}.
\end{thm}
\begin{proof}
For the linear operator $\Pi_h$ in Lemma \ref{verf1} we choose a Scott-Zhang operator. Then Lemma \ref{verf1} and Corollary \ref{verf4} easily imply the following discrete inf-sup condition
\begin{equation} \label{LBBbarVbarQ}
  \sup_{0 \neq v_h \in \bar{V}_h}
  \frac{ b(v_h,q_h)}{\|v_h\|_1} \gtrsim \|q_h\|_0 
        \quad \text{for all} \ q_h \in \bar{Q}_h. 
\end{equation}
Then the discrete LBB condition easily follows for the case (a) of pure Dirichlet boundary conditions and case (b) of mixed boundary conditions, since $V_h = \bar{V}_h$ and $Q_h \subset \bar{Q}_h$ in these cases. 

For the case (c) of pure Neumann boundary conditions it follows from \eqref{LBBbarVbarQ} that, for each $q_h \in \bar{Q}_h = Q_h$, there exists a $\bar{v}_h \in \bar{V}_h$ such that
\[
  b(\bar{v}_h,r_h) = (q_h,r_h) \quad \text{for all} \ r_h \in Q_h \quad \text{and} \quad \|\bar{v}_h\|_1 \lesssim \|q_h\|_0.
\]
Let $c \in P_0(\Omega)^d$ be the mean value of $\bar{v}_h$ on $\Omega$ and set $v_h = \bar{v}_h - c$. Then $v_h \in V_h$ and
\[
  b(v_h,r_h) = b(\bar{v}_h,r_h) = (q_h,r_h) \quad \text{for all} \ r_h \in Q_h \quad \text{and} \quad \|v_h\|_1 \le \|\bar{v}_h\|_1 \lesssim \|q_h\|_0,
\]
which immediately implies the discrete LBB condition in this case.
\end{proof}

\subsection{The discrete BP condition}\label{subsec2}

Next we will discuss inf-sup conditions for a different pair of norms. First we show a local version of \cite[Proposition 1]{BercovierPironneau_1979}.

\begin{thm} \label{BerPirlocal}
{Let $k \in \mathbb{N}$ with $k \ge 2$, let $K \in \mathcal{T}_h$ and assume that Condition \eqref{restriction} holds for $K$.}
Then
\[
  \sup_{0 \neq v \in \bar{V}_K}
  \frac{ b_K(v,q)}{\|v\|_{0,K}} \gtrsim \|\grad q\|_{0,K}
  \quad \text{for all} \ q \in \bar{Q}_K.
\]
\end{thm}
\begin{proof}
By using Lemma \ref{Wienerslocal} we obtain
\[
   \sup_{0 \neq v \in \bar{V}_K}
  \frac{ b_K(v,q)}{\|v\|_{0,K}}
  \ge \frac{ b_K(T_K q,q)}{\|T_K q\|_{0,K}}
  \gtrsim \frac{h_K^2 \,\|\grad q\|_{0,K}^2}{h_K^2 \, \|\grad q\|_{0,K}}
    = \|\grad q\|_{0,K}.
\]
\end{proof}

The next result is a global version of Theorem \ref{BerPirlocal} and  can be seen as an extension of \cite[Proposition 1]{BercovierPironneau_1979}, where the proof was given for $k=2$ and for meshes made of rectangles for $d=2$ and bricks for $d=3$.

\begin{thm} \label{BerPirglobal}
{Let $k \in \mathbb{N}$ with $k \ge 2$ and assume that Condition \eqref{restriction} holds for all $K \in \mathcal{T}_h$.}
Then
\[
  \sup_{0 \neq v_h \in V_h}
  \frac{ b(v_h, q_h)}{\|v_h\|_0} \gtrsim \|\grad q_h\|_0 
        \quad \text{for all} \ q_h \in Q_h.
\]
\end{thm}
\begin{proof}
We will restrict ourselves to the case $d=3$, the case $d=2$ is a simple adaptation of the arguments used for $d=3$ and is omitted.

In a first step we show the discrete BP condition in spaces without scaling conditions, i.e.,
\begin{equation} \label{barVbarQ}
  \sup_{0 \neq v_h \in \bar{V}_h}
  \frac{ b(v_h, q_h)}{\|v_h\|_0} \gtrsim \|\grad q_h\|_0 
        \quad \text{for all} \ q_h \in \bar{Q}_h .
\end{equation}
Let $\mathcal{N}_h$ be the set of all nodes of the mesh $\mathcal{T}_h$, and
let $v^{(\node)}$ and $\varphi_{\node}$ be the nodal value and the nodal basis function associated to a node $\node \in \mathcal{N}_h$, respectively.
Then
\[
  \|v_h\|_0^2 
  \sim \sum_{\node \in \mathcal{N}_h} 
    \|v^{(\node)} \varphi_{\node}\|_0^2
  \quad \text{for} \ v_h
    =  \sum_{\node \in \mathcal{N}_h} 
    v^{(\node)} \, \varphi_{\node}.
\]
For $q_h \in \bar{Q}_h$ it follows that
\begin{align*}
  \sup_{0 \neq v_h \in \bar{V}_h}
  \frac{b(v_h,q_h)^2}{\|v_h\|_0^2}
   \sim \sup_{0 \neq v_h \in \bar{V}_h}
        \frac{\left[ \sum_{\node \in \mathcal{N}_h} b(v^{(\node)} \, \varphi_{\node},q_h)\right]^2
             }%
             {\sum_{\node \in \mathcal{N}_h} \|v^{(\node)} \,  \varphi_{\node}\|_0^2
             } 
  = \sum_{\node \in \mathcal{N}_h}
    \sup_{0 \neq v^{(\node)} \in \mathbb{R}^3}
  \frac{b(v^{(\node)} \varphi_{\node},q_h)^2}{\|v^{(\node)} \varphi_{\node}\|_0^2}.
\end{align*}
We discuss three cases for the node $\node$.
\begin{enumerate}
\item
Assume that $\node \in \mathcal{N}_\text{int}(K)$ for some $K \in \mathcal{T}_h$.
By setting 
\[
  v^{(\node)} = \left( J \grad \hat{q}_h\right)(\noderef)
\]
we obtain as in the proof of Lemma \ref{Wienerslocal}
\[
   b(v^{(\node)} \varphi_{\node},q_h)  
    = b_K(v^{(\node)} \varphi_{\node},q_h)  
    \sim \, h_K^3 \, \|\grad \hat{q}_h(\noderef)\|^2
\]
and
\[
  \|v^{(\node)} \varphi_{\node}\|_0^2 
   = \|v^{(\node)} \varphi_{\node}\|_{0,K}^2 
   \sim h_K^3 \, \|v^{(\node)}\|^2 
    \lesssim h_K^{5} \, \|\grad \hat{q}_h(\noderef)\|^2,
\]
which implies
\[
  \sup_{0 \neq v^{(\node)} \in \mathbb{R}^3}
  \frac{b(v^{(\node)} \varphi_{\node},q_h)^2}{\|v^{(\node)} \varphi_{\node}\|_0^2}
  \gtrsim h_K \, \|\grad \hat{q}_h(\noderef)\|^2 .
\]
\item
Next assume that $\node \in \mathcal{N}_\text{int}^*(\mathcal{E}_i)$ for some $K \in \mathcal{T}_h$. 
By setting
\[
  v^{(\node)} = \left(  J_i \, \partial_i \hat{q}_h 
        \right)(\noderef) 
\]
we obtain as in the proof of Lemma \ref{Wienerslocal}
\[
   b_K(v^{(\node)} \varphi_{\node},q_h)  
    \sim \, h_K^3 \, |\partial_i \hat{q}_h(\noderef)|^2
\]
and
\[
  \|v^{(\node)} \varphi_{\node}\|_{0,K}^2 
    \lesssim h_K^{5} \, |\partial_i \hat{q}_h(\noderef)|^2.
\]
Analogous estimates follow for all other elements $K'$ with $\node \in K'$. This leads to 
\[
   b(v^{(\node)} \varphi_{\node},q_h)
     = \sum_{K \in \mathcal{T}(\node)} b_K(v^{(\node)} \varphi_{\node},q_h)
    \sim \sum_{K \in \mathcal{T}(\node)} h_K^3 \, |\partial_i \hat{q}_h(\noderef)|^2
\]
and
\[
  \|v^{(\node)} \varphi_{\node}\|_0^2 
   = \sum_{K \in \mathcal{T}(\node)} \|v^{(\node)} \varphi_{\node}\|_{0,K}^2 
    \lesssim \sum_{K \in \mathcal{T}(\node)} h_K^5 \, |\partial_i \hat{q}_h(\noderef)|^2,
\]
with $\mathcal{T}(\node) = \{K \in \mathcal{T}_h \colon \node \in K\}$.
It easily follows that
\[
  \sup_{0 \neq v^{(\node)} \in \mathbb{R}^3}
  \frac{b(v^{(\node)} \varphi_{\node},q_h)^2}{\|v^{(\node)} \varphi_{\node}^2\|_0}
  \gtrsim \sum_{K \in \mathcal{T}(\node)} h_K \, |\partial_i \hat{q}_h(\noderef)|^2 .
\]
\item
Finally, assume that $\node \in \mathcal{N}_\text{int}^*(\mathcal{F}_{i,j})$ for some $K \in \mathcal{T}_h$. Then it follows completely analogously to the previous case that
\[
  \sup_{0 \neq v^{(\node)} \in \mathbb{R}^3}
  \frac{b(v^{(\node)} \varphi_{\node},q_h)^2}{\|v^{(\node)} \varphi_{\node}\|_0^2}
  \gtrsim \sum_{K \in \mathcal{T}(\node)} 
    h_K \, \left( |\partial_i \hat{q}_h(\noderef)|^2 + |\partial_j \hat{q}_h(\noderef)|^2 \right) .
\]
\end{enumerate}
In summary, we have
\begin{align*}
  \sup_{0 \neq v_h \in \bar{V}_h} \frac{b(v_h,q_h)^2}{\|v_h\|_0^2} 
  & \gtrsim \sum_{K \in \mathcal{T}_h}
      h_K \, \bigg( \sum_{\node \in \mathcal{N}_\text{int}(K)}
        \|\grad \hat{q}_h(\noderef)\|^2  + 
      \sum_{\substack{\node \in \mathcal{N}_\text{int}^*(\mathcal{F}_{i,j})\\i,j \in \{1,2,3\}, i < j}}
        \left(|\partial_i \hat{q}(\noderef)|^2 + |\partial_j \hat{q}(\noderef)|^2 \right)
        \nonumber \\  
  & \quad \qquad\qquad\qquad {} + \sum_{\substack{\node \in \mathcal{N}_\text{int}^*(\mathcal{E}_i)\\i \in \{1,2,3\}   }}
        |\partial_i \hat{q}(\noderef)|^2  \bigg) 
      \sim \sum_{K \in \mathcal{T}_h} 
         \|\grad q_h\|_{0,K}^2 =  \|\grad q_h\|_0^2 ,
\end{align*}
which completes the proof of \eqref{barVbarQ}. 

Then the discrete BP condition easily follows for the case (a) of pure Dirichlet boundary conditions and the case (b) of mixed boundary conditions, since $V_h = \bar{V}_h$ and $Q_h \subset \bar{Q}_h$ in these cases. The discrete BP condition for the case (c) of pure Neumann boundary conditions can be shown completely analogously as for the discrete LBB condition in Theorem \ref{discreteLBB}. 
\end{proof}

\section{Conclusions and open problems}\label{sec4}

Previous work on the discrete LBB condition of the Stokes problem for the generalized Taylor-Hood family $Q_k$--$Q_{k-1}$ on quadrilateral/hexahedral meshes focused on the case of pure Dirichlet boundary conditions. In this paper we extended the analysis to the case, where "do-nothing" boundary conditions are prescribed on part of the boundary or on the whole boundary.
This was one of two aims of the paper. A second aim was the analysis of the discrete BP condition, which plays an important role in the analysis of the smoothing property of some multigrid method for the discretized Stokes problem. 

From the technical point of view the major difference to previous work is the element-wise approach as opposed to the widely used macroelement technique. The element-wise approach allowed to derive element-wise discrete LBB and BP conditions.

The present paper suffers from the same rather severe restrictions on hexahedral meshes in 3D as in previous work. The analysis of discrete inf-sup conditions for general hexahedral meshes remains an open problem. Another open problem is the analysis of isoparametric generalized Taylor-Hood families in 2D and 3D to cope with curved boundaries. Perturbation arguments similar to those used in \cite{Bressan_2011}, \cite{Bressan_2013} for isogeometric generalized Taylor-Hood families seem to be a promising approach for this open problem.

\begin{appendices}

\section{The continuous LBB condition}\label{secA1}
The LBB condition \eqref{classicalinfsup} can be rewritten as
\[
  \|q\|_0 \le c \, \|\grad_V \, q\|_{V^*} \quad \text{for all} \ q \in Q
\]
with $c = 1/\beta$, the linear functional
\[
  \langle \grad_V \, q,v \rangle := - (q,\div v)_0 \quad \text{for} \ v \in V,
\]
and the dual norm
\[
  \|\grad_V \, q\|_{V^*} = \sup_{0 \neq v \in V}
                  \frac{\langle \grad_V \, q,v \rangle}{\|v\|_1}.
\]
If $V = H_0^1(\Omega)^d$ the linear functional $\grad_V \, q$ can be identified with the gradient of $q$ in the distributional sense. In this case we omit the subscript $V$ and simply write $\grad q$.

We have the following important inequality:
\begin{lem}[Ne\v{c}as]
Let $\Omega \subset \mathbb{R}^d$ be a bounded and open set with Lipschitz continuous boundary. Then there exists a constant $c_N > 0$, such that 
\begin{equation} \label{infsup4}
  \|q\|_0 \le c_N \, (\|q\|_{-1} + \|\grad q\|_{-1} \quad \text{for all} \ 
  q \in L^2(\Omega).
\end{equation}  
\end{lem}

\begin{proof}
See \cite{Necas12}, under stronger assumptions also \cite{Duvaut76}.
\end{proof}

The LBB condition for $V = H_0^1(\Omega)^d$ and $Q = L_0^2(\Omega)$ follows by standard compactness arguments.

\begin{thm}
Let $\Omega \subset \mathbb{R}^d$ be a bounded, connected and open set with Lipschitz continuous boundary. Then there exists a constant $c > 0$, such that
\begin{equation} \label{infsup5}
  \|q\|_0 \le c \, \|\grad q\|_{-1} \quad \text{for all} \ q \in L_0^2(\Omega). 
\end{equation}  
\end{thm}
For later reference we recall the well-known proof.
\begin{proof}
Assume the inequality (\ref{infsup5}) is not valid. Then there exists a sequence $(q_k)$ in $L_0^2(\Omega)$ with $\|q_k\|_0 = 1$ and
$\|\grad q_k\|_{-1} \to 0$. 
Since the embedding $i \colon H_0^1(\Omega) \longrightarrow L^2(\Omega)$ is compact, the (adjoint) embedding $i^* \colon L^2(\Omega) \longrightarrow H^{-1}(\Omega)$ is also compact.
Therefore, there exists a convergent sub-sequence $(q_k')$ in
$H^{-1}(\Omega)$. From (\ref{infsup4}) it follows that 
\begin{align*}
  \|q_k' - q_\ell'\|_0 
    & \le c_N \, (\|q_k' - q_\ell'\|_{-1} + \|\grad (q_k' - q_\ell')\|_{-1}) \\
    & \le c_N \, (\|q_k' - q_\ell'\|_{-1} + \|\grad q_k'\|_{-1} + \|\grad q_\ell'\|_{-1})
      \to 0,
\end{align*}
if $k, \ell \to \infty$. So $(q_k')$ is a Cauchy sequence in $L^2(\Omega)$ and, therefore, $q_k' \to q$ in $L^2(\Omega)$ for some $q \in  L_0^2(\Omega)$.

For $v \in \left(H_0^1(\Omega)\right)^d$ we have:
\begin{align*}
  |(q,\div v)_0| = \lim_{k\to \infty} |(q_k',\div v)_0| 
   = \lim_{k\to \infty} |\langle \grad q_k',v \rangle| 
   \le \lim_{k\to \infty} \|\grad q_k'\|_{-1} \, \|v\|_1 = 0.
\end{align*}
Hence 
\[
  (q,\div v)_0 = 0 \quad \text{for all} \ v \in \left(H_0^1(\Omega)\right)^d,
\]
i.e., the generalized gradient of $q$ vanishes, which implies that $q$ is constant. 
Since $q \in L_0^2(\Omega)$ it follows $q = 0$, in contradiction to $\|q\|_0 = \|q_k\|_0 = 1$. 
\end{proof}

The LBB condition can be shown in a similar way 
for the spaces $V = H_{0,\Gamma_D}^1(\Omega)^d$ and $Q = L^2(\Omega)$ in the case $|\Gamma_N| > 0$:

\begin{thm} \label{GammaD}
Let $\Omega \subset \mathbb{R}^d$ be a bounded, connected and open set with Lipschitz continuous boundary, and let $V = H_{0,\Gamma_D}^1(\Omega)^d$ with $|\Gamma_N| > 0$. Then there exists a constant $c > 0$, such that
\begin{equation} \label{infsup7}
  \|q\|_0 \le c \, \|\grad_V \, q\|_{V^*} \quad \text{for all} \ q \in L^2(\Omega). 
\end{equation}  
\end{thm}
\begin{proof}
Assume the inequality (\ref{infsup7}) is not valid. Then there exists a sequence $(q_k)$ in $L^2(\Omega)$ with $\|q_k\|_0 = 1$ and $\|\grad_V \, q_k\|_{V^*} \to 0$. Since $\|\grad q_k\|_{-1} \le \|\grad_V \, q_k\|_{V^*}$, it follows that $\|\grad q_k\|_{-1} \to 0$.
As in the proof of Theorem it follows that
there exists a sub-sequence $(q_k')$ with $q_k' \to q$ in $L^2(\Omega)$ for some $q \in  L^2(\Omega)$.

For $v \in V$ we have:
\begin{align*}
  |(q,\div v)_0| = \lim_{k\to \infty} |(q_k',\div v)_0| 
   = \lim_{k\to \infty} |\langle \grad_V \, q_k',v \rangle| 
   \le \lim_{k\to \infty} \|\grad_V \, q_k'\|_V^* \, \|v\|_1 = 0.
\end{align*}
Hence 
\[
  (q,\div v)_0 = 0 \quad \text{for all} \ v \in V.
\]
Since $H_0^1(\Omega)^d \subset V$, the generalized gradient of $q$ vanishes, which implies that $q$ is constant. 
Moreover, it additionally follows from  above
\[
  (q,\div v)_0 = q \, \int_\Omega \div v \ d x = 0 \quad \text{for all} \ v \in V.
\]
Since there exists a $v \in V$ with $\int_\Omega \div v \ d x \neq 0$, see next lemma, it follows that $q = 0$, in contradiction to $\|q\|_0 = \|q_k\|_0 = 1$. 
\end{proof}

\begin{lem}
Under the assumptions of Theorem \ref{GammaD} there exists a function $v \in V$ such that
\[
  \int_\Omega \div v \ d x \neq 0.
\]
\end{lem}
\begin{proof}
According to Lemma 1.5.1.9 in \cite{grisvard11} there exists $w \in C^\infty(\overline\Omega)$ such that
\[
  w \cdot n \ge \delta  \quad \text{almost everywhere on} \ \Gamma
\]
for some $\delta > 0$. Let $\rho \in C^\infty(\mathbb{R}^d)$ be a non-negative mollifier with center at some point on $\Gamma_N$ and with a sufficiently small support. Then $v = \rho \, w \in V$ and
\[
  \int_\Omega \div v \ d x = \int_\Gamma \rho \, w \cdot n \ d \sigma \ge \delta \int_\Gamma \rho \ d \sigma > 0. 
\]
\end{proof}

An alternative proof of Theorem \ref{GammaD} can be found in \cite{Baerland_2017}. The proof there also requires the result of the last lemma, which was not explicitly addressed in \cite{Baerland_2017}.

Finally, for the case of pure Neumann boundary conditions we have
\begin{thm} \label{contLBBpureNeumann}
Let $\Omega \subset \mathbb{R}^d$ be a bounded, connected and open set with Lipschitz continuous boundary, and let $V = \hat{H}^1(\Omega)^d$. Then there exists a constant $c > 0$, such that
\[
  \|q\|_0 \le c \, \|\grad_V \,  q\|_{V^*} \quad \text{for all} \ q \in L^2(\Omega). 
\]  
\end{thm}
The LBB condition for this case of pure Neumann boundary conditions can be shown completely analogously as for the discrete LBB condition, see Theorem \ref{discreteLBB}, and is omitted. 

\section{Some useful relations}\label{secA2} 

For the convenience of the reader we recall several useful relations, which are used throughout the paper.

\subsubsection*{Norms under the change of variables}

For all $q \in H^1(K)$:
\[
  \|\grad q\|_{0,K}^2 
     \sim h_K^{d-2} \, \|\grad \hat{q}\|_{0,\hat K}^2,\quad
  \|q\|_{0,f}^2 
     \sim h_K^{d-1} \, \|\hat{q}\|_{0,\hat f}
\]
For all $v \in H^1(K)^d$:
\begin{align*}
  \|v\|_{0,K}^2 
     \sim h_K^d \, \|\hat{v}\|_{0,\hat K}^2, \quad
  \|v\|_{1,K}^2 
     \lesssim h_K^{d-2} \, \|\hat{v}\|_{1,\hat K}^2
\end{align*}

\subsubsection*{Norms in terms of nodal values}

For all $\hat{q} \in Q_{k-1}(\hat K)$:
\begin{align*}
  \|\hat{q}\|_{0,\hat K}^2 
     \sim \sum_{\noderef \in \mathcal{N}(\hat K)} |\hat{q}(\noderef)|^2, \quad
  \|\grad \hat{q}\|_{0,\hat K}^2
     \sim \sum_{\noderef \in \mathcal{N}(\hat K)} \|\grad \hat{q}(\noderef)\|^2, \quad
  \|\hat{q}\|_{0,\hat f}^2
     \sim \sum_{\noderef \in \mathcal{N}(\hat f)} |\hat{q}(\noderef)|^2
\end{align*}
For all $\hat{v} \in Q_k(\hat K)^d$:
\begin{align*}
  \|\hat{v}\|_{0,\hat K}^2
     \sim \sum_{\noderef \in \mathcal{N}(\hat K)} \|\hat{v}(\noderef)\|^2, \quad
  \|\hat{v}\|_{1,\hat K}^2
     \sim \sum_{\noderef \in \mathcal{N}(\hat K)} \|\hat{v}(\noderef)\|^2
\end{align*}

\subsubsection*{Relations based on Gauss-Lobatto quadrature}

For all $\hat g \in Q_{2k-1}(\hat K)$ with $\hat g(\noderef) \ge 0$ for all $\noderef \in \mathcal{N}(\hat K)$:
\begin{align*}
  \int_{\hat K} \hat g \ d \hat x
    & \sim \sum_{\noderef \in \mathcal{N}(\hat K)} \hat g(\noderef)
\end{align*}
For all $\hat g \in Q_{2k-1}(\hat f)$  with $\hat g(\noderef) \ge 0$ for all $\noderef \in \mathcal{N}(\hat f)$:
\begin{align*}
  \int_{\hat f} \hat g \ d \hat \sigma
    & \sim \sum_{\noderef \in \mathcal{N}(\hat f)} \hat g(\noderef)
\end{align*}
For all $\hat g \in Q_{2k-1}(\hat f)$:
\begin{align*}
  \Big|\int_{\hat f} \hat g \ d \hat \sigma \Big|
    & \lesssim \sum_{\noderef \in \mathcal{N}(\hat f)} |\hat g(\noderef)|
\end{align*}
\end{appendices}

\subsubsection*{Acknowledgements:}

The author is very grateful for helpful information on the topic of the paper
provided by Andrea Bressan, Vivette Girault (communicated by Thomas Wick), Rolf Rannacher (communicated by Thomas Wick), David Silvester, Rolf Stenberg, and Christian Wieners.

\bibliographystyle{abbrvurl}
\bibliography{infsupTaylorHood}

\end{document}